\documentclass[twocolumn,a4paper,aps,prl,amssymb,amsmath,floatfix]{revtex4}

\newcommand{\Z}{\mathbb{Z}}
\newcommand{\Q}{\mathbb{Q}}
\newcommand{\F}{\mathbb{F}}
\newcommand{\T}{T}
\newcommand{\ra}{\longrightarrow}

\newcommand{\op}{{\oplus}}

\newtheorem{theorem}{Theorem}
\begin{document}

\hyphenation{co-kernel}

\title{Integer Cech Cohomology of Icosahedral Projection Tilings}

\author{F.~G\"ahler\footnote{Present address: Fakult\"at f\"ur Mathematik, 
  Universit\"at Bielefeld, Postfach 10 01 31, D-33615 Bielefeld, Germany
  \newline Email: gaehler@math.uni-bielefeld.de
}}
\affiliation{
  Institut f\"ur Theoretische und Angewandte Physik, 
  Universit\"at Stuttgart, D-70550 Stuttgart, Germany
}
\author{J.~Hunton}
\affiliation{
  Department of Mathematics, University of Leicester,
  Leicester, LE1 7RH, United Kingdom 
}
\author{J.~Kellendonk}
\affiliation{
  Universit\'e de Lyon, Universit\'e Lyon 1, 
  CNRS, UMR5208, Institut Camille Jordan,
  43 blvd du 11 novembre 1918, F-69622 Villeurbanne cedex, France
}

\begin{abstract}
  The integer Cech cohomology of canonical projection tilings of dimension 
  three and codimension three is derived. These formulae are then evaluated
  for several icosahedral tilings known from the literature. Rather 
  surprisingly, the cohomologies of all these tilings turn out to have 
  torsion. This is the case even for the Danzer tiling, which is, in some
  sense, the simplest of all icosahedral tilings. This result is in contrast
  to the case of two-dimensional canonical projection tilings, where many 
  examples without torsion are known.
\end{abstract}

\maketitle

\section{Introduction}
\label{sec:intro}


The local isomorphism class (LI class) of a repetitive tiling can be given a
metric, which makes it a compact topological space of tilings, called the hull
of the tiling.  Many important properties of the tiling depend on the topology
of its hull, which can be characterized to some extent by topological
invariants. For a brief introduction to the topology of tiling spaces, and why
it matters, we refer to \cite{Sad06}. The simplest topological invariants are
the Cech cohomology groups, which moreover are practically computable for
large classes of tilings. This relies mainly on the fact that tiling spaces
are inverse limits of simpler spaces \cite{Sad03}. Using such an inverse limit
construction, Anderson and Putnam \cite{AP} gave a practical algorithm to
compute the Cech cohomology of primitive substitution tilings, an algorithm
which works for integer, rational, or any other kind of coefficients.  Another
important class of tilings whose cohomology is computable are canonical 
projection tilings, with a definition of `canonical' \cite{FHKcmp} which is 
actually more general than the usual one. Basically, one only has to assume 
polyhedral windows whose faces satisfy certain rationality conditions. For 
such canonical projection tilings, the cohomology with rational coefficients 
is determined by the formalisms of Forrest, Hunton and Kellendonk (FHK)
\cite{FHKcmp,FHKmem} and of Kalugin \cite{Kalugin}. For the co-dimension 2
case, even the integer result was implicitly contained in \cite{FHKcmp}, but
the torsion part remained unnoticed because it was thought
\cite{FHKcmp,FHKmem} that the cohomology of canonical projection tilings was
always free. This turned out to be incorrect, however, and was disproved by
counterexamples \cite{FGBanff}.  Whereas the co-dimension 2 case is therefore
easy to extend to integer coefficients (also in Kalugin's formalism), the
co-dimension 3 case had resisted any such attempt for a long time. Recent
progress enables us, however, to overcome these difficulties in both
formalisms. The refined theories of FHK and Kalugin both allow now to
determine the cohomology with integer coefficents for a large number of
examples, even though there remain cases where the torsion part cannot be
fully determined. In this paper, we sketch the essential steps and ingredients
of the computation for an integer version of Kalugin's formalism, whose setup
is somewhat easier to explain. The results for a number of well-known
icosahedral tilings from the literature are then presented. It turns out that
all these tilings have torsion in their integer Cech cohomology. For a more
detailed description, and for an integer version of the FHK formalism, we
refer to a forthcoming publication \cite{GHK}.


\section{Canonical Projection Tilings}
\label{sec:ProTil}

The tilings we consider here are known as canonical projection tilings, 
where we use the broader, more general notion of `canonical' from 
\cite{FHKcmp}, rather than the traditional one used in \cite{FHKmem}. 
Despite their name, canonical projections tilings are more conveniently
described as {\em sections} through a higher-dimensional periodic tiling. 
For this purpose, we consider a Euclidean space $E^{n}=E^{\parallel} \oplus
E^{\perp}$, containing a lattice $\Gamma$ of full rank $n$, whose projection
on either summand is both bijective and dense. In $E^{n}$, we have a
$\Gamma$-periodic tiling, whose tiles $K_{i}$ are prisms
$K_{i}=K_{i}^{\parallel} \times K_{i}^{\perp}$.  Generic sections parallel to
$E^{\parallel}$ yield then tilings of $E^{\parallel}$ by tiles
$K_{i}^{\parallel}$. Each tile $K_{i}^{\parallel}$ has an acceptance domain
$K_{i}^{\perp}\subset E^{\perp}$: whenever the section touches it, the tile is
present in the tiling. Sections touching the boundary of an acceptance domain
are called singular; they correspond to more than one tiling. The periodic
tiling can be wrapped on the $n$-torus $\T^{n}=E^{n}/\Gamma$. Disregarding
singular sections, the points in this torus label the possible tilings in our
tiling space $\Omega$. This is known as the torus parametrization
\cite{BHP97}.

In the following, we consider tilings whose acceptance domains are
polyhedral. If a section touches a face $F$ in the boundary of an 
acceptance domain $K^{\perp}$, the corresponding tiling is not unique
in the set $F \times K^{\parallel}$ (which has codimension 1). Often, 
more than one tile is singular at a time. The ambiguity of a singular 
section has to be resolved by cutting up the torus along the singular 
set where the tiling is ambiguous, and completing it from both sides.
However, a point in the cut-up torus still does not necessarily 
correspond to a non-singular cut. The tiling may still be singular 
somewhere away from the origin. To avoid this, one would have to cut up 
the torus along the whole hyperplane $F \times E^{\parallel}$. Kalugin 
\cite{Kalugin} has proposed to proceed in steps, using an inverse limit 
construction: at each step, the torus is cut up only by a finite amount.
One thereby arrives at a sequence of tori $\Omega_{r}=\T^{n}\setminus A_{r}$, 
cut-up along larger and larger sets $A_{r}$. Kalugin has shown that the 
full tiling space $\Omega$ is given by the inverse limit of this sequence 
of cut-up tori $\Omega_{r}$. Moreover, their cohomology satisfies the
exact sequence
\begin{multline}
  \xrightarrow{} H^{k}(\Omega_{r}) 
  \xrightarrow{} H_{n-k-1}(A_{r}) \\ \xrightarrow{} 
  H_{n-k-1}(\T^{n}) \xrightarrow{} H^{k+1}(\Omega_{r}) \xrightarrow{} 
  \label{eq:ex1}
\end{multline}
Next, we make use of the fact that canonical projection tilings have
acceptance domains with rationally oriented faces (with respect to the
projected lattice $\Gamma^{\perp}$). For the icosahedral tilings 
considered here, for which $n=6$, this implies that any component 
$F \times K^{\parallel}$ of the singular set $A_{r}$ is contained in a 
$5$-dimensional hyperplane comprising 4 independent lattice directions. 
This set is therefore 4-periodic, containing whole 4-dimensional 
arrays of sets $F \times K^{\parallel}$, and if $r$ is increased, these 
pieces grow together to form, for some finite $r$, a full, rational 
4-plane thickened along a fifth direction. Actually, this thickened 
rational 4-plane is wrapped onto the torus, and the set $A_{r}$ is 
therefore a finite union of rational 4-tori, each thickened along a 
fifth direction. The homotopy type of $A_{r}$ (and thus its homology) 
does not change if $r$ is further increased. Moreover, it is often 
possible to find a finite arrangement $\tilde A$ of normal, thin 4-tori,
which is homotopy equivalent to $A_{r}$, and thus has the same homology 
(see \cite{Kalugin}). For all the 3-dimensional icosahedral tilings 
considered below, this is the case, and we can work with such a set 
$\tilde A$, which consists of a union of $L_{2}$ 4-tori $\T^{\alpha}$, 
indexed by $\alpha\in I_{2}$, each of which is stabilized by the 
translations in a sublattice $\Gamma^{\alpha}$ of rank 4. These 4-tori 
intersect in $L_{1}$ 2-tori $T^{\theta}$, indexed by $\theta\in I_{1}$, 
each stabilized by a sublattice $\Gamma^{\theta}$ of rank 2, and 
$L_{0}$ intersection points. The $L_{1}^{\alpha}$ 2-tori contained in 
$\T^{\alpha}$ are indexed by $I_{1}^{\alpha}$, the $L_{0}^{\alpha}$ 
intersection points in $\T^{\alpha}$ by $I_{0}^{\alpha}$, and the 
$L_{0}^{\theta}$ intersection points in $\T^{\theta}$ by $I_{0}^{\theta}$. 
Even though the set $\tilde A$ to be removed from the torus $\T^{6}$ is 
already reasonably simple, we have to make one more replacement by yet 
another homotopy equivalent set: we pass from $\tilde A$ to its simplicial
resolution $A$ (see \cite{Kalugin}). The set $A$ is the union of a still 
larger number of spaces, which comprises also spaces homotopic to the 
former intersection tori. The advantage of this is, that the maximal number 
of spaces with a common intersection is drastically reduced, which proves 
useful in the computation of the homology of $A$ via a Mayer-Vietoris 
spectral sequence.


\section{Cohomology Computation}
\label{sec:Cohom}

\begin{table*}
\begin{center}
{
\renewcommand{\arraystretch}{1.4}
\begin{tabular}{|c|c|c|}
\hline
$\oplus_{\alpha\in I_{2}}\Lambda_{4}\Gamma^{\alpha}$ & & \\
\hline
$\oplus_{\alpha\in I_{2}}\Lambda_{3}\Gamma^{\alpha}$ & & \\
\hline
$\oplus_{\alpha\in I_{2}}\Lambda_{2}\Gamma^{\alpha}\bigoplus
 \oplus_{\theta\in I_{1}}\Lambda_{2}\Gamma^{\theta}$
  & $\oplus_{\alpha\in I_{2}}\oplus_{\theta\in I_{1}^{\alpha}}\Lambda_{2}\Gamma^{\theta}$ & \\
\hline
$\oplus_{\alpha\in I_{2}}\Lambda_{1}\Gamma^{\alpha}\bigoplus
 \oplus_{\theta\in I_{1}}\Lambda_{1}\Gamma^{\theta}$
  & $\oplus_{\alpha\in I_{2}}\oplus_{\theta\in I_{1}^{\alpha}}\Lambda_{1}\Gamma^{\theta}$ & \\
\hline
$\Z^{L_{2}}\bigoplus \Z^{L_{1}} \bigoplus \Z^{L_{0}}$ & 
$\oplus_{\alpha\in I_{2}}\Z^{L_{1}^{\alpha}} \bigoplus
 \oplus_{\alpha\in I_{2}}\Z^{L_{0}^{\alpha}} \bigoplus 
 \oplus_{\theta\in I_{1}}\Z^{L_{0}^{\theta}} $ &
$\oplus_{\alpha\in I_{2}} \oplus_{\theta\in I_{1}^{\alpha}}\Z^{L_{0}^{\theta}}$ \\
\hline
\end{tabular}
}
\end{center}
\caption{First page of the Mayer-Vietoris spectral sequence for the homology 
         of $A$.}
\label{tab:MVSS}
\end{table*}

We have now seen that the singular set $A_{r}$ is, for sufficiently large
$r$, homotopy equivalent to the set $A$, which can therefore replace $A_{r}$
in the exact sequence (\ref{eq:ex1}), and after the disappearance of $r$
the direct limit of the cohomology of $\Omega_{r}$ is easily carried out.
Moreover, the sequence (\ref{eq:ex1}) can be split as follows:
\begin{multline}
0 \xrightarrow{} S_{k}(R) \xrightarrow{} H^{k}(\Omega;R) 
\xrightarrow{} H_{6-k-1}(A;R) \\ \xrightarrow{\alpha^{k+1}_{R}} 
  H_{6-k-1}(\T^{6};R) \xrightarrow{} S_{k+1}(R) \xrightarrow{} 0
\label{eq:ex2}
\end{multline}
where we have explicitly indicated the coefficient group $R$. To determine
$H^{k}(\Omega;R)$, we need to compute $H_{*}(\T^{6};R)$, $H_{*}(A;R)$, and
$S_{k}(R) = \text{coker}\, \alpha^{k}_{R}$ for the different coefficient 
groups $R=\Q$, $\Z$, or $\F_{p}$. A few conclusions on $S_{k}(R)$ can 
immediately be drawn. Since $H^{k}(\Omega;R)=0$ for $k>3$, we have 
$S_{k}(R)=0$ for $k>3$, and $\alpha^{4}_{R}$ is onto. Also, because $A$ 
has dimension $4$, we get $S_{k}(R)=H_{6-k}(\T^{6};R)$ for $k<2$, and 
since the constituent spaces of $A$ have intersections of dimension at
most $2$, $S_{2}(R)$ is also easily computed: $S_{2}(R) = \Lambda_{4}
\Gamma\otimes R / \langle \Lambda_{4}\Gamma^{\alpha}\otimes R \rangle$,
where $\Lambda_{k}\Gamma$ denotes the $k$-fold totally antisymmetric tensor
product of $\Gamma$.

Before we address the computation of $S_{3}(R)$, we determine the homology 
of $A$. As $A$ is a finite union of spaces $A_{i}$, its homology can be
computed using a Mayer-Vietoris spectral sequence (MVSS), which expresses 
the homology of $A$ in terms of the homology of the $A_{i}$ and their 
intersections. We use a homological version of the MVSS used in \cite{Basu}, 
where the reader can also find more information on the evaluation of the 
MVSS in a similar setting. It is important here that we have passed to the 
simplicial resolution $A$ of $\tilde{A}$, which drastically reduces the 
number of columns required in the MVSS (which is equal to the maximal
number of spaces $A_{i}$ having a common intersection, compare 
\cite{Kalugin}). The first page (see \cite{Basu}) $E^{1}_{k,\ell}$ of 
the Mayer-Vietoris double complex contains the direct sums of the 
$\ell$-homologies of all intersections of $k$ constituent spaces $A_{i}$
of $A$. It is given in Table~\ref{tab:MVSS}. All these groups are free, 
and we also know that $H_{0}(A)=R$ and $H_{1}(A)=R^{6}$. As a consequence, 
the total rank of all differentials in this and the second page, departing 
from the diagonal $k+\ell=2$, can be computed. Beyond this total rank, 
explicit knowledge of those differentials is not needed. Only one 
differential needs to be known explicitly: 
\[
d^{1}_{{1,2}}: 
  \underset{\alpha\in I_{2}}{\bigoplus} \oplus_{\theta\in I_{1}^{\alpha}}
                \Lambda_{2} \Gamma^{\theta} \rightarrow 
(\oplus_{\alpha\in I_{2}^{\ }} \Lambda_{2} \Gamma^{\alpha})
\,\oplus\,  
(\oplus_{\theta\in I_{1}^{\ }} \Lambda_{2} \Gamma^{\theta}) \notag
\]
Its source space is embedded in the obvious way in both sectors
of the target space, with different signs. The cokernel of the differential 
$d^{1}_{{1,2}}$ is the only possible source of torsion in $H_{*}(A;\Z)$.
As it appears in the leftmost column of the double complex, the homology
of $A$ is still the direct sum of $E^{\infty}_{k,\ell}$, $n=k+\ell$,
despite the presence of torsion. The final result on the homology of $A$ 
is summarized in Theorem~\ref{thm1}.


\begin{theorem}
Suppose $A$ is an arrangement of singular tori as above. For any 
coefficient group $R=\Q$, $\Z$, or $\F_{p}$, we then have:
\begin{align}
H_{0}(A;R) &= R \notag\\
H_{1}(A;R) &= R^{6} \notag \\
H_{2}(A;R) &= \textrm{coker}\  d^{1}_{1,2}(R) \oplus R^{f} \notag \\ 
H_{3}(A;R) &= (\oplus_{\alpha\in I_{2}}\Lambda_{3}\Gamma^{\alpha}) \oplus 
              (\textrm{ker}\ d^{1}_{1,2}(R)) \notag \\
H_{4}(A;R) &= \oplus_{\alpha\in I_{2}}\Lambda_{4}\Gamma^{\alpha} \notag \\
\text{with}\quad f &= -3 L_{2} - L_{1} + \sum_{\alpha\in I_{2}} L_{1}^{\alpha}
                      + 5 + \chi 
     \notag\\
\text{and}\quad \chi &= L_{0} - \sum_{\alpha\in I_{2}} L_{0}^{\alpha}
     + \sum_{\alpha\in I_{2}} \sum_{\theta\in I_{1}^{\alpha}} L_{0}^{\theta}
     - \sum_{\theta\in I_{1}} L_{0}^{\theta} \notag
\end{align}
For $R=\Z$, the only potential source of torsion 
in $H_*(A;\Z)$ is $\textrm{coker}\ d^{1}_{1,2}(\Z)$.
$\chi$ is the Euler characteristic of $\Omega$.
\label{thm1}
\end{theorem}


Next, we need to determine $S_{3}(R)$ for the different coefficient groups
$R$ (all other $S_{k}(R)$ are already known). We start with $R=\Q$. For 
$S_{3}(\Q)$ we need the image of $\alpha_{\Q}^{3}$, which is at least 
$M = \langle \Lambda_{3}\Gamma^{\alpha}\otimes\Q \rangle $. We claim
that it cannot be bigger than $M$. To see this, consider a simplicial 
decomposition of $\T^{6}$, containing $A$ as a subcomplex, and take the 
chain maps
\[
C_{3}(A) \overset{i}{\ra} C_{3}(\T^{6}) \overset{P}{\ra} 
\Lambda_{3}\Gamma\otimes\Q
\]
$P$ assigns to each simplex its volume form. As $P$ vanishes on boundaries,
it induces an isomorphism $H_{3}(\T^{6})\rightarrow\Lambda_{3}\Gamma\otimes\Q$.
All cells in $A$ map under $P$ to directions contained in $M$, and so the 
image of $\alpha_{\Q}^{3}$ cannot be larger than $M$. The rational result 
gives us bounds on the possible values of $\text{Im}\,\alpha_{\Z}^{3}$:
\[
M_{1}= \langle \Lambda_{3}\Gamma^{\alpha}\rangle 
\subseteq \text{Im}\,\alpha_{\Z}^{3}
\subseteq \langle \Lambda_{3}\Gamma^{\alpha}\otimes\Q\rangle 
\cap \Lambda_{3}\Gamma = M_{2} \notag
\]
If $M_{1}=M_{2}$, $S_{3}(\Z) = \text{coker}\, \alpha_{\Z}^{3} = 
\Lambda_{3}\Gamma / \langle \Lambda_{3}\Gamma^{\alpha}\rangle $ is free,
and we are done. Otherwise, we can try the sharper lower bound
\begin{equation}
M_{1}' = \langle b_{1}^{\theta}\wedge b_{2}^{\theta}\wedge v | 
  \theta \in I_{1}, \langle b_{1}^{\theta}, b_{2}^{\theta}\rangle =
  \Gamma^{\theta}, v\in\Gamma\rangle
\label{eq:bound2}
\end{equation}
on $\text{Im}\,\alpha_{\Z}^{3}$, which comes about as follows. The 
$3$-chains in $C_{3}(A)$ need not be closed within a single constituent 
torus of $A$. Rather, the boundaries of $3$-chains contained in different 
$4$-tori may cancel each other. To see this, we unfold the torus for a 
moment, and consider two lattice $4$-planes $p_{1}$ and $p_{2}$, intersecting
in a $2$-plane $\ell_{1}$ with stabilizer $\Gamma^{\theta}$. Let $v$ be any 
lattice vector. The planes $p_{1}$ and $p_{2}+v$ intersect in a plane 
$\ell_{2}$ parallel to $\ell_{1}$, and the planes $p_{1}+v$ and $p_{2}+v$ 
intersect in a plane $\ell_{3}=\ell_{1}+v$. We can now construct two 
$3$-chains, one contained in $p_{1}$ with boundary $[\ell_{2}]-[\ell_{1}]$,
and one contained in $p_{2}$ with boudnary $[\ell_{3}]-[\ell_{1}]$. As
$\ell_{1}$ and $\ell_{3}$ are equivalent in the torus, the boundaries
of the two $3$-chains cancel, and we have a $3$-cyle homologous to 
$b_{1}^{\theta}\wedge b_{2}^{\theta}\wedge v$, where $\{b_{i}^{\theta}\}$ 
is a basis of $\Gamma^{\theta}$; this proves the bound (\ref{eq:bound2}).
Again, if $M'_{1}=M_{2}$, $S_{3}(\Z) = \text{coker}\, \alpha_{\Z}^{3} = 
\Lambda_{3}\Gamma / M_{2} $ is free.


We now restrict ourselves to the case where $S_{3}(\Z)$ is free.
There may be still torsion in $H_{2}(A;\Z)$, however, and it is not
clear yet whether it lifts to $H^{3}(\Omega;\Z)$. To see this, note 
that the sequence $0 \rightarrow \Z \rightarrow X \rightarrow \Z_{2} 
\rightarrow 0$ has two solutions, $X=\Z$ and $X=\Z\oplus \Z_{2}$.
No integral arguments will help to distinguish these cases, 
and we need further input. It turns out that the knowledge of
$S_{3}(\F_{p})$ for suitable primes $p$ will help. Our aim is to 
show that if $S_{3}(\Z)$ is free, then $rk(S_{3}(\Q)) = rk(S_{3}(\F_{p}))$. 
To do so, we write $A \xrightarrow{} \T^{6}$ as a composite:
\[
A \xrightarrow{\Delta} A\times \ldots \times A 
\xrightarrow{a_{1}\times \ldots \times a_{6}}
S^{1}\times\ldots\times S^{1} \xrightarrow{i} \T^{6}
\]
In homology, we then get
\begin{multline}
H_{3}(A) \xrightarrow{\Delta*} \oplus_{r_{1}+\cdots +r_{6}=3}\, 
H_{r_{1}}(A)\otimes \cdots \otimes H_{r_{6}}(A) \\
\xrightarrow{a_{1}*\otimes \cdots \otimes a_{6}*} 
\oplus_{r_{1}+\cdots +r_{6}=3}\, 
H_{r_{1}}(S^{1})\otimes\cdots\otimes H_{r_{6}}(S^{1}) \\
\xrightarrow{i*} H_{3}(\T^{6}) \notag
\end{multline} 
Since $H_{k}(S^{1})=0$ for $k\ge2$, only $H_{0}(A)$ and $H_{1}(A)$ can
contribute. For those dimensions, however, $H_{k}(A)$ is free; passing 
from $\Q$ to $\F_{p}$ can therefore not change the rank of 
$\text{Im}\,\alpha^{3}$.


\begin{table*}
\begin{center}
{
\renewcommand{\arraystretch}{1.4}
\begin{tabular}{|c|c|c|c|c|c|c|l|}
\hline
$H^3$                          & $H^2$                      & $H^1$     & $H^0$ & $\chi$ & planes   & $\Gamma$ &  \\
\hline
$\Z^{  20}\op\Z_2$             & $\Z^{ 16}$                 & $\Z^{ 7}$ & $\Z$  &     10 & 5-fold   & F        & Danzer  \\
\hline
$\Z^{ 181}\op\Z_2$             & $\Z^{ 72}\op\Z_2$          & $\Z^{12}$ & $\Z$  &    120 & mirror   & P        & Ammann-Kramer \\
\hline
$\Z^{ 331}\op\Z_2^{20}\op\Z_4$ & $\Z^{102}\op\Z_2^4\op\Z_4$ & $\Z^{12}$ & $\Z$  &    240 & mirror   & F        & dual canonical $D_6$ \\
\hline
$\Z^{ 205}\op\Z_2^{2}$         & $\Z^{ 72}$                 & $\Z^{ 7}$ & $\Z$  &    145 & 3,5-fold & F        & canonical $D_6$  \\
\hline
\end{tabular}
}
\end{center}
\caption{Cohomology of icosahedral tilings from the literature.
  Also indicated are the Euler characteristic $\chi$, the type of 
  icosahedral lattice (P for primitive, F for face-centered), 
  and the orientation of the boundaries of the acceptance domains 
  (perpendicular to 3-fold or 5-fold axes of the icosahedron, 
  or parallel to its mirror planes).}
\label{tab:ex}
\end{table*}

We are now going to compare the ranks of the groups in (\ref{eq:ex2}) 
for the different coefficient groups $\Q$ and $\F_{p}$. From the universal 
coefficient theorem we know that
\begin{align}
rk(H^{k}(&\Omega;\F_{p})) - rk(H^{k}(\Omega;\Q))\notag\\ &= 
T_{p}(H^{k}(\Omega;\Z)) + T_{p}(H^{k+1}(\Omega;\Z)),  \notag\\
rk(H_{k}(&A;\F_{p})) - rk(H_{k}(A;\Q))\notag\\ &=
T_{p}(H_{k}(A;\Z)) + T_{p}(H_{k-1}(A;\Z)), \notag
\end{align}
where $T_{p}(X)$ is the rank of $p$-torsion in $X$. Clearly, we are
interested mainly in those primes $p$ for which $H_{2}(A;\Z)$ has 
$p$-torsion. Inserting this into (\ref{eq:ex2}), and using 
$rk(S_{3}(\Q)) = rk(S_{3}(\F_{p}))$, we obtain
\begin{gather}
 rk(H^{3}(\Omega;\F_{p})) - rk(H^{3}(\Omega;\Q)) = T_{p}(H_{2}(A;\Z)),\notag\\ 
 \text{Torsion}(H^{3}(\Omega;\Z)) = \text{Torsion}(H_{2}(A;\Z)). \notag
\end{gather}
We therefore have arrived at our main result, which is summarized in 
Theorem~\ref{thm2}:


\begin{theorem}
Consider a dimension 3, co-dimension 3 canonical projection tiling
with arrangement of singular tori $A$ as above, and suppose $S_{3}(\Z)$ 
is free. For any coefficient field $R=\Q$ or $\F_{p}$, we then have
\begin{multline}
rk(H^{k}(\Omega;R)) = rk(H_{6-k-1}(A;\Z)) + rk(S_{k}(R)) \\
 + rk(S_{k+1}(R)) - \tbinom{6}{k+1}. \notag 
\end{multline}
For the coefficient ring $R=\Z$, the same relation holds for the free ranks,
and  
\begin{align}
 \text{Torsion}(H^{2}(\Omega;\Z)) &= \text{Torsion}(S_{2}(\Z))   \notag \\
 \text{Torsion}(H^{3}(\Omega;\Z)) &= \text{Torsion}(H_{2}(A;\Z)) \notag
\end{align}
All other $H^{k}(\Omega;\Z))$ are free.
\label{thm2}
\end{theorem}


\section{Results and Discussion}
\label{sec:ResDisc}

We have evaluated the formulae from Theorems 1 and 2 for several 
icosahedral tilings from the literature, with the help of a computer.
These tilings comprise the Danzer tiling \cite{Dan89}, the Ammann-Kramer 
tiling \cite{KN84}, the canonical $D_{6}$ and the dual canonical $D_{6}$ 
\cite{KP95}. The results are given in Table~\ref{tab:ex}. Most remarkably, 
all these tilings have torsion in $H^{3}(\Omega)$, whose origin is in the 
torsion of $H_{2}(A)$. Two examples have torsion also in $H^{2}(\Omega)$, 
resulting from the torsion in $S_{2}(\Z)$. In two dimensions, at least 
the simpler canonical projection tilings have no torsion, so that it is 
rather surprising that even the Danzer tiling with its combinatorially 
simple set $A$ has torsion. It should be noted, however, that the torsion 
of the Danzer tiling has a different origin. It comes from the homology of 
$A$, which does not provide any torsion in the two-dimensional case. 
It should also be remarked that $S_{3}$ is free for all these examples. 
Only for the dual canonical $D_{6}$ we had to resort to the sharper bound 
(\ref{eq:bound2}) to show this.

The methods presented here are applicable also in the 2-dimensional case,
but everything is considerably simpler. Details on the 2-dimensional case
will be presented elsewhere \cite{GHK}.


\end{document}